# Rejoinder

**Donald B. Rubin**

It is a pleasure to have the opportunity to respond to these three complementary discussions written by Edward L. Korn (ELK), Paul R. Rosenbaum (PRR) and Stephen E. Fienberg (SEF), each of whom has made substantial contributions to problems of causal inference. Thanks to all three for the generosity expressed in their comments. I am in agreement with many points raised by them, but there are a few places where we may not fully agree, possibly due to misunderstandings.

## ELK

ELK's ordering of problems by their importance, displayed in his Figure 1, is difficult to dispute. I like to think of there being a continuum in causal inference from perfect randomized experiments to sloppy observational studies, and I like to "chip away" at all points in the continuum because I often feel that there is potentially a tremendous transfer of insights from work on one problem to work on another. For a specific example, the work on noncompliance presented in Angrist, Imbens and Rubin (1996) led to the general principal stratification framework in Frangakis and Rubin (2002), which provided a formal structure for embedding the resolution of the "censoring/truncation due to death" problem used a decade earlier in the context of an actual FDA submission, as noted in Rubin (1998, 2000) and the target article.

ELK's examples are highly appropriate and offer strong support for the importance of working in the context of real problems. I particularly liked his orthodontic example from Korn, Teeter and Baumrind (2001), and think that generalizations of the idea used there could be quite fruitful, especially because that idea implies specific suggestions for improvements to the design of particular types of observational studies. ELK is correct that when the principal strata are effectively observed, as in that example, distributional assumptions, such as normality, can be avoided.

As pointed out in Angrist, Imbens and Rubin (1996) and its rejoinder, the focus of our analysis of noncompliance is on the units in this experiment who would comply when assigned either treatment. This group is, by definition, the only collection of units in this experiment who can be observed receiving and not receiving the treatment, and thus is the only collection of units providing any data about the causal effect of receiving versus not receiving the treatment. The always-takers are always observed receiving the treatment, and the never-takers are never observed receiving the treatment, and thus data from neither of these latter two groups provides any evidence about the effect of taking versus not taking the treatment. In other words, in this experiment "efficacy" can only be estimated for the compliers. "Effectiveness" is a joint property of (a) efficacy for the compliers, (b) how the treatment is "marketed" (i.e., how compliance is enforced) and (c) "placebo" effects of assignment on the noncompliers (i.e., the always-takers and the never-takers). This point is discussed at length in Sheiner and Rubin (1995).

Consequently, the emphasis within the principal stratification framework is on separately estimating efficacy and the other components of effectiveness, in order to help the generalization to other situations with possibly different placebo effects or different marketing effects. Thus I am puzzled by ELK's statement that the principal stratification approach "... negates one of the usual reasons for being interested in efficacy and not effectiveness..."—quite the opposite in my mind, but perhaps I misunderstood his meaning.

## PRR

I do understand ELK's desire to avoid the entire problem of censoring due to death in the quality-of-life example by assigning the lowest possible QOL score to those who are dead, but as I have argued, this approach, to me, mixes up issues of (a) estimating the scientific effect of a treatment intervention in groups of units where it can estimated and (b) individual value judgments about the value of death versus various qualities of life. PRR seems to agree, even providing a Seneca quote in support! If we accept ELK's suggestion to use, for example, a ranksum test, it seems an approach such as the one advocated by PRR is quite attractive because it avoids having to make a particular choice of a single value for the QOL of someone who is dead, common for everyone, and thus PRR's approach allows diverging patient preference orderings of death relative to different qualities of life.

As is typical with PRR's contributions, I find his proposal deep and creative. I do wonder, however, about its implications for applied consumers of such data, such as doctors or patients. PRR's proposal does avoid the need for distributional assumptions, noted by ELK and presented in my Section 6, but the proposal seems to replace them with confidence intervals for estimands that may not be as easily understood, that is, for the order statistics that would have been observed if all $n$ treated subjects had instead received control. Despite the undeniable mathematical elegance of the approach, and the clear exposition conveyed especially by the subsequent examples, my sympathies continue to be with what I view as the more direct Bayesian model-based formulation. But this preference may be largely a matter of taste and differing experiences.

Also, to some extent I realize that this view about PRR's proposal being difficult to convey to consumers may be unfair because the unfamiliar is nearly always seen as more challenging, and I look forward to seeing a family of work evolving from this interesting idea. For example, what would this approach have to say about the simpler noncompliance problem, in particular when we impose the no-defier assumption, but we do not impose the exclusion restrictions for both always-takers and never-takers, as in Hirano et al. (2000)?

## SEF

SEF points out that my article largely avoided the use of formal notation and equations, in contrast to some of my earlier work. I agree with SEF that such formality is the best way to nail down intuitive ideas, but I think that we also agree that sometimes informal exposition works better for conveying the underlying ideas to less technical audiences. The lectures on which this paper was based were delivered in the late afternoon to large audiences, which included nonstatistical relatives of Morrie and Morris, and I did not want to put them to sleep! I hope that, if the words in the written version are not precise enough to convey critical ideas clearly, the technical references given will make up for any deficiencies.

Regarding terminology, for years I have avoided the use of "counterfactuals" to describe "potential outcomes" for two major reasons. First, at the design phase, no well-defined potential outcome is counterfactual, although there do exist a priori counterfactuals, such as the value of an outcome when exposed to treatment for a never-taker who will never be exposed to the treatment no matter the assignments. Second, at the analysis phase, at least some of the potential outcomes are factual. I agree with SEF's preference for making all potential outcomes random variables, as I did in Rubin (1975, 1978). But to bridge the Bayesian potential outcomes framework to non-Bayesians, it is important to recognize that in Neyman's original formulation, which is the classical randomization-based formulation still used by many, if not most, statisticians today, for example by PRR, the potential outcomes are not random variables; instead the potential outcomes are treated as "...fixed features of the finite population of $N$ subjects." In this approach, only the randomization indicator is a random variable.

SEF's plea at the conclusion of his Section 2 argues for formality of causal inference using graphical models and, as stated earlier, I certainly agree that having more formality available is better (as in the full Bayesian approach to noncompliance developed in Imbens and Rubin, 1997). However, as pointed out in Rubin (2004)—including in my rejoinder to Lauritzen (2004), for causal inference I find the graphical approach more ambiguous and more confusing, as well as less flexible and less formal, than the potential outcomes approach. But maybe this too is simply a matter of differing tastes.

An aspect of SEF's discussion that I find puzzling is his comment about the "average causal effect (ACE)," and his implication that in Rubin (1978) I focused on the ACE and justified randomization to a Bayesian using the ACE. But that *Annals* article was devoted to discussing the full posterior predictive distribution of all the potential outcomes, from which all Bayesian causal inferences follow, no matter how the causal effects are defined; for example, see (4.1) in that article. And in Rubin (1974, pages 690–694) I pointed out that the use of the average causal effect was rather arbitrary but convenient for frequentist justifications for randomization based on unbiased estimation over the randomization distribution, and I noted the possible use of the median or mid-mean instead (page 690). Over the years, many researchers (e.g., Brillinger, Jones and Tukey, 1978) have used definitions of causal effects other than the "ACE," a term coined, I believe, in Holland (1986). Also see the definitions of causal effects on quantiles implied by PRR's confidence intervals in his equation (1) here, which clearly do not equal the ACE.

## CONCLUSION

In conclusion, I once again thank all three discussants for their comments and hope that this package assembled by the editorial board of *Statistical Science* will lead to more work in this great research area, much of it following up on the interesting ideas contributed by the discussants.